\documentclass{article}
\usepackage[]{epsf,epsfig,amsmath,amssymb,amsfonts,latexsym}
\title{Double-Tail Invariant Measures of the Dyck Shift}
\author{Tom Meyerovitch}

\newtheorem{theorem}{Theorem}[section]
\newtheorem{lemma}{Lemma}[section]
\newtheorem{proposition}{Proposition}[section]

\newenvironment{proof}{{\bf Proof:} \rm}{\hfill $\Box$ \medskip\\}

\begin{document}
\maketitle
\begin{abstract}
In \cite{T04} it was shown that the One sided Dyck is uniquely
ergodic with respect to the one sided-tail relation, where the
tail invariant probability is also shift invariant and obtains the
topological entropy.\\
In this paper we show that the two sided Dyck has a double-tail
invariant probability, which is also shift invariant, with entropy
strictly less than the topological entropy.
\end{abstract}

\section{Introduction}
Let $\Sigma$ be a finite alphabet. For a subshift $X \subset
\Sigma^\mathbb{Z}$, we define the double-tail relation, or
\emph{homoclinic} \cite{PS97} relation of $X$ to be:
\[\mathcal{T}_2(X) := \{(x,x') \in X \times X\:\ \exists n \geq 0 \ \forall |k|>n \ x_k=x'_k \}\]
A $\mathcal{T}_2(X)$-holonomy is an injective Borel function $g:A
\mapsto g(A)$, with $A$ a Borel set and $(x,g(x)) \in
\mathcal{T}_2(X)$ for every $x \in A$. We say that $\mu \in
\mathcal{M}(X)$ is a double-tail invariant if $\mu(A)=\mu(g(A))$
for every $\mathcal{T}_2(X)$-holonomy $g$.\\
In this paper we characterize the double-tail invariant
probabilities for the Dyck shifts. In addition to its two
equilibrium measures, the two sided Dyck shift has another
double-tail invariant probability -- shift invariant,
non-equilibrium. These are the only three double-tail invariant,
ergodic probabilities on the two sided Dyck shift.\footnote{This
article is a part of the author's M.Sc. thesis, written under the
supervision of J. Aaronson, Tel-Aviv University}

\subsection{Definition of the Dyck Shift}
Fix an integer $m>1$. Throughout the rest of this paper we denote:
\[\Sigma=\{\alpha_1,\ldots,\alpha_m\}\cup\{\beta_1,\ldots,\beta_m\}\]
The \emph{Dyck monoid} $M$, is the monoid (with $0$), with
generators $\alpha_i, \beta_i$, $1\leq i \leq m$, and relations:
\begin{enumerate}
\label{def_m}
  \item \(\alpha_j \cdot \beta_j \equiv \Lambda \equiv 1 (mod M) , j=1, \ldots ,m\)
  \item \(\alpha_i \cdot \beta_j \equiv 0 (mod M) ,i \neq j \)
\end{enumerate}

Note that by the very definition of $M$ as a monoid (a semigroup
with natural element), if $w\equiv w'$ (mod $M$), then for every
$s,t \in \Sigma^*$, $swt \equiv sw't$ (mod $M$).

The \textsl{m-Dyck Language} is
\[ L=\{l \in \Sigma^* : l \neq 0 (mod M) \} \]
and the corresponding \text{(two sided) m-Dyck subshift} is
\[ X=\{ x \in \Sigma^\mathbb{Z} : (x_i)_{i=r}^{l} \in L \mbox{ for
all } -\infty  <r \leq l < +\infty \} \]
These are indeed subshifts, since we only pose restrictions on finite blocks.\\
Note that the 1-Dyck shift is simply the full 2-Shift, and so we
will only be interested in the case where $m\geq 2$.\\

 For $x \in X$, let
\[ H_i(x)= \left\{ \begin{array}{ll}
  \sum_{j=0}^{i-1}
\sum_{l=1}^{m}(\delta_{\alpha_l,x_j}-\delta_{\beta_l,x_j}) &
\mbox{if
$i>0$} \\
\sum_{j=i}^{-1}\sum_{l=1}^{m}(\delta_{\beta_l,x_j}-\delta_{\alpha_l,x_j})
& \mbox{if $i<0$} \\
0 & \mbox{if $i=0$} \\
\end{array}
\right.
\]

\section{Maximal Entropy Implies Double-Tail Invariance}\label{mu_+_mu_-}
In \cite{WK74} it was demonstrated that the Dyck shift has two
ergodic shift invariant probabilities with entropy equal to the
topological entropy. Such probabilities are called
\emph{equilibrium states}. In this section we show that both of
these probabilities are also shift invariant.\\
 We introduce the following sets, which are mutually disjoint and are
tail-invariant. For $s,t \in \{\{+\infty\},\{-\infty\}, \mathbb
R\}$ we define:
\[B^s_t=\{x \in X : \liminf_{i \rightarrow + \infty}H_i(x) \in s ,
    \liminf_{i \rightarrow - \infty}H_i(x) \in t \}\]
let
\[\Omega^+_-=\{x \in \{\alpha_1 \ldots \alpha_m,\beta\}^{\mathbb{Z}} : \liminf_{i \rightarrow + \infty}\widehat{H}_i(x)=+\infty ,
 \liminf_{i \rightarrow - \infty}\widehat{H}_i(x)=-\infty \}\]
 and
 \[\Theta^-_+=\{x \in \{\beta_1 \ldots \beta_m,\alpha\}^{\mathbb{Z}} : \liminf_{i \rightarrow + \infty}\widetilde{H}_i(x)=-\infty ,
 \liminf_{i \rightarrow - \infty}\widetilde{H}_i(x)=+\infty \}\]
 Where $\widehat{H}$ and $\widetilde{H}$ are the cocycles
 generated by $1_{\bigcup[\alpha_j]}-1_{[\beta]}$ and
 $1_{[\alpha]}-1_{\bigcup[\beta_j]}$ respectively.
 Define:
  \[g_+:B^+_- \mapsto \Omega^+_-\]
 \[ (g_+(y))_{i}= \left\{
\begin{array}{ll}
\alpha_j & y_{i}=\alpha_j  \\
\beta & y_{i} \in \{ \beta_{1},\ldots,\beta_{m}\}
\end{array}
\right.
\]
\[g_-:B^-_+ \mapsto \Theta^-_+\]
 \[ (g_-(y))_{i}= \left\{
\begin{array}{ll}
\beta_j & y_{i}=\beta_j  \\
\alpha & y_{i} \in \{ \alpha_{1},\ldots,\alpha_{m}\}
\end{array}
\right.
\]
$g_+$ is a Borel bijection from $B^{+\infty}_{-\infty}$ to
$\Omega^{+\infty}_{-\infty}$ and $g_-$ is a Borel bijection of the
appropriate sets. The definitions of $g_+$ and $g_-$ can also be
extended to functions $g_+:B^{\mathbb{R}}_{-\infty} \mapsto
\Omega^{\mathbb{R}}_{-\infty}$ and $g_-:B^{-\infty}_{\mathbb{R}}
\mapsto \Theta^{-\infty}_{\mathbb{R}}$, which are also Borel
bijections.

\begin{lemma}\label{g_+_g_-_isomorphism}
$g_+:B^{+\infty}_{-\infty} \mapsto \Omega^{+\infty}_{-\infty}$,
$g_-:B^{-\infty}_{+\infty} \mapsto \Theta^{-\infty}_{+\infty}$,
$g_+:B^{\mathbb{R}}_{-\infty} \mapsto
\Omega^{\mathbb{R}}_{-\infty}$, $g_-:B^{-\infty}_{\mathbb{R}}
\mapsto \Theta^{-\infty}_{\mathbb{R}}$ are isomorphisms of the two
sided tail relations:
\[ (g_+ \times g_+) (\mathcal{T}_2(B^{+\infty}_{-\infty}))=
\mathcal{T}_2(\Omega^{+\infty}_{-\infty})\]
\[ (g_- \times g_-) (\mathcal{T}_2(B^{-\infty}_{+\infty}))=
\mathcal{T}_2(\Theta^{-\infty}_{+\infty})\]
\[ (g_+ \times g_+) (\mathcal{T}_2(B^{\mathbb{R}}_{-\infty}))=
\mathcal{T}_2(\Omega^{\mathbb{R}}_{-\infty})\]
\[ (g_- \times g_-) (\mathcal{T}_2(B^{-\infty}_{\mathbb{R}}))=
\mathcal{T}_2(\Theta^{-\infty}_{\mathbb{R}})\]

\end{lemma}

\begin{proof}
We prove the result for $g_+:B^{+\infty}_{-\infty} \mapsto
\Omega^{+\infty}_{-\infty}$, the other results are proved in the same manner.\\
$ (g_+ \times g_+) (\mathcal{T}_2(B^{+\infty}_{-\infty})) \subset
\mathcal{T}_2(\Omega^{+\infty}_{-\infty})$ is trivial, so we show
the other inclusion. Suppose $(g_+(x),g_+(y)) \in
\mathcal{T}_2(\Omega^{+\infty}_{-\infty})$. Let $n_0 \geq 0$ be
such that $g_+(x)_{[-n_0,n_0]^c}=g_+(y)_{[-n_0,n_0]^c}$.\\
Let \[r(i,x)=max\{j<i: H_j(x)=H_i(x)\}\] Clearly,
$r(i_1,x)=r(i_2,x)$ is impossible
 for $i_1 \neq i_2$.
 Since
$\liminf_{n\rightarrow +\infty}H_n(x)>-\infty$ and $\liminf_{n
\rightarrow +\infty}H_n(y)>-\infty$, there exists $c$ such that
for some large $N$, $H_i(x)>c$ for every $i>N$. Since $\liminf_{n
\rightarrow -\infty}H_n(x) = \liminf_{n \rightarrow
-\infty}H_n(y)= -\infty$ , it follows that there exist some
$i_0<N$ such that $H_{i_0}(x)=c$, so for every $i>N$,
$r(i,x)>i_0$. The same argument applies for $y$. Since
$(r(i,x))_{i>N}$ and $((r(i,y))_{i>N}$ are both injective
sequences of integers, bounded from below, it follows that
\[\lim_{n \rightarrow +\infty}r(n,x)=\lim_{n \rightarrow
+\infty}r(n,y)=+\infty\]
 Note that for $n_1,n_2> n_0$,
\[\widehat{H}_{n_1}(g(x))-\widehat{H}_{n_2}((g(x))=\widehat{H}_{n_1}(g(y))-\widehat{H}_{n_2}((g(y))\]
 so for all large $n$ enough so that $r(n,x)>n_0$,$r(n,y)>n_0$, either:
 \begin{enumerate}
    \item $g(x)_n=g(y)_n=m+1$,in which case $r(n,x)=r(n,y)$ and $x_{r(n,x)}=y_{r(n,y)}$, so $x_n=y_n$
    \item $g(x)=g(y)=i$ for $1< i < m$, and then $x_n= y_n=\alpha_i$
 \end{enumerate}
 Obviously, for $n< -n_0$, $x_n=y_n$. This proves $(x,y) \in
 \mathcal{T}_2(B^{+\infty}_{-\infty})$.
\end{proof}

\begin{lemma}\label{mu_+_mu_-_unique}
There exists a unique $\mathcal{T}_2$-invariant probability of $X$
supported by $B^{+\infty}_{-\infty}$, and a unique
$\mathcal{T}_2$-invariant probability of $X$ supported by
$B^{-\infty}_{+\infty}$. There are no $\mathcal{T}_2$-invariant
probabilities on $B^{\mathbb{R}}_{-\infty}$ and
$B^{-\infty}_{\mathbb{R}}$.
\end{lemma}

\begin{proof}
The symmetric product measure $p$ on $\Omega$ assigns probability
one to $\Omega^+_-$ Transporting the product measure on $\Omega$
by means of $g^{-1}$ to $B^+_-$ yields a tail invariant
probability measure on $X$, by the previous lemma.\\
On the other hand, any tail invariant probability on $X$ supported
by $B^+_-$ can be transported to a tail invariant probability on
$\Omega$ by $g$, and by the uniqueness of tail invariant
probability on $\Omega$, we conclude the uniqueness of double tail
invariant probability on $B^+_-$. This also proves that no double
tail invariant probability on  $B^{\mathbb{R}}_{-\infty}$ exist.
We obtain the results for $B^-_+$ and $B^{-\infty}_{\mathbb{R}}$
symmetrically.
\end{proof}

\section{A Third Double-tail Invariant
Probability}\label{third_prob}
 For $z \in \{0,1\}^\mathbb{Z}$, we define:
\[ \tilde{H}_i(z)= \left\{ \begin{array}{ll}
  \sum_{j=0}^{i-1}
(\delta_{1,z_j}-\delta_{0,z_j}) & \mbox{if
$i>0$} \\
\sum_{j=i}^{-1}(\delta_{0,z_j}-\delta_{1,z_j})
& \mbox{if $i<0$} \\
0 & \mbox{if $i=0$} \\
\end{array}
\right.
\]
 Let \[S_{-\infty}^{-\infty}= \{z \in \{0,1\}^\mathbb{Z} : \;
\liminf_{n \rightarrow + \infty}\tilde{H}_n(z)=-\infty \;
\liminf_{n \rightarrow -\infty}\tilde{H}_n(z)=-\infty \;\}\] Let
us define a Borel function $F: S_{-\infty}^{-\infty} \times
\{1,\ldots,m\}^\mathbb{Z}
\mapsto \Sigma^\mathbb{Z}$:\\

Let
\[F(z,a)_n= \left\{
\begin{array}{ll} \alpha_j & \mbox{if $z_n=1$, $a_{\gamma_n(z)}=j$} \\
\beta_j & \mbox{if $z_n=0$, $k=\varepsilon_n(z)$, and
$a_{\gamma_k(z)}=j$}
\end{array}
\right.
\]
where,
\[\gamma_k(z) = \left\{ \begin{array}{ll}
 \sum_{i=0}^{k}z_i & k \geq 0\\
 -\sum_{i=k}^{-1}z_i  & k<0\\
\end{array}
\right.\]
\[\varepsilon_n(z) = \max\{l < n : \;
\tilde{H}_l(z) \leq\tilde{H}_{n+1}(z)\}\]
 Since $\liminf_{n \rightarrow -\infty}\tilde{H}_n(z)=-\infty$
for $z \in S_{-\infty}^{-\infty}$, $F$ is well defined.\\
\begin{lemma}
For every $z \in S_{-\infty}^{-\infty}$, $a \in
\{1,\ldots,m\}^\mathbb{Z}$, $F(z,a) \in X$.
\end{lemma}
\begin{proof}
Suppose $x=F(z,a) \not\in X$, then there exist $n,n' \in
\mathbb{Z}$, $n<n'$, such that $x_n=\alpha_i$, $x_{n'}=\beta_j$
with $i \neq j$ and $n=\max\{l < n' : \; H_l(x)=H_{n'+1}(x)\}$.
But in that case, $n=\varepsilon_{n'}(z)$, so
$i=j=a_{\gamma_n(z)}$.
\end{proof}
 Let $\mu_1$ be the symmetric product measure on $\{0,1\}^\mathbb{Z}$, and
$\mu_2$ the symmetric product measure on
$\{1,\ldots,m\}^\mathbb{Z}$.
\begin{lemma}
$\mu_1(S_{-\infty}^{-\infty})=1$
\end{lemma}
\begin{proof}
This follows from the ergodicity of the skew-product
$\{0,1\}^{\mathbb{Z}} \times \mathbb{Z}$ given by the cocycle
$\tilde{H}_k$ with respect to the product measure $\mu_1 \times
\nu_{\mathbb{Z}}$ where $\nu_{\mathbb{Z}}$ is the counting measure
on $\mathbb{Z}$.
\end{proof}
We define: $\tilde{\mu}=(\mu_1 \times \mu_2)\circ F^{-1}$. Since
$F^{-1}(B_{-\infty}^{-\infty})=S^{-\infty}_{-\infty} \times
\{1,\ldots,m\}^\mathbb{Z}$ it follows
that $\tilde{\mu}(B_{-\infty}^{-\infty})=1$.\\

Let us also define a Borel mapping $z:B_{-\infty}^{-\infty}
\mapsto S_{-\infty}^{-\infty}$:
\[z(x)_n= \left\{ \begin{array}{ll} 1 & x_n \in \{\alpha_1,\ldots,\alpha_m\} \\
0 & x_n \in \{\beta_1,\ldots,\beta_m\} \end{array} \right. \]

The following lemma gives an explicit formula the $\tilde{\mu}$
probability of a cylinder:
\begin{lemma}\label{tilde_mu_explicit}
Let $w \in L(X)$. If the number of paired $\alpha$'s in $w$ is
$n_1$ and the number of unpaired $\alpha$'s and $\beta$'s is $n_2$
($2n_1+n_2=|w|$) then
$\tilde{\mu}([w]_k)=m^{-(n_1+n_2)}(\frac{1}{2})^{|w|}$.
\end{lemma}
\begin{proof}
Denote by $f_1,\ldots,f_{n_1}$ the locations of matched $\alpha$'s
in $w$. Denote by $g_1,\ldots,g_{n'_2}$ the locations of unmatched
$\alpha$'s in $w$. Denote by $h_1,\ldots,h_{n''_2}$ the locations
of unmatched $\beta$'s in $w$. We have $n'_2+n''_2=n_2$. For
$\vec{r} \in \mathbb{Z}^{n_1}$, $\vec{s} \in \mathbb{Z}^{n'_2}$,
$\vec{t} \in \mathbb{Z}^{n''_2}$, define:
\[A_{\vec{r}}=\{z :\; \gamma_{k+f_l}(z)=r_l \; 1\leq l\leq n_1 \}\]
\[B_{\vec{s}}=\{z :\; \gamma_{k+g_l}(z)=s_l \; 1\leq l\leq n'_2 \}\]
\[C_{\vec{t}}=\{z :\; \gamma_{\varepsilon_l}(z)=t_l \; \varepsilon_l=\varepsilon_{k+h_l}(z) 1\leq l\leq n''_2 \}\]
Informally, $A_{\vec{r}},B_{\vec{s}},C_{\vec{t}}$ determine the
locations in the sequence $a \in \{1,\ldots,m\}^{\mathbb{Z}}$
involved in selecting the types of $\alpha$'s and $\beta$'s within
the coordinates $k,\ldots,k+|w|$. Now we define:
\[Z = \{z \in S_{-\infty}^{-\infty} :\; z_{i+k} = z(w)_i \mbox{ for } 0 \leq i
\leq |w| \}\]
\[A'_{\vec{r}}=\{a \in \{1,\ldots,m\}^{\mathbb{Z}}:\; a_{r_l}=j \mbox{ if } w_{f_l}=\alpha_j\}\]
\[B'_{\vec{s}}=\{a \in \{1,\ldots,m\}^{\mathbb{Z}}:\; a_{s_l}=j \mbox{ if } w_{g_l}=\alpha_j\}\]
\[C'_{\vec{t}}=\{a \in \{1,\ldots,m\}^{\mathbb{Z}}:\; a_{t_l}=j \mbox{ if } w_{h_l}=\beta_j\}\]
With the above definitions we can write:
\begin{equation}
F^{-1}([w]_k)= Z \times \{1,\ldots,m\}^{\mathbb{Z}}\cap
\bigcup_{\vec{s},\vec{t},\vec{r}}((A_{\vec{r}}\times
A'_{\vec{r}})\cap(B_{\vec{s}}\times B'_{\vec{s}}) \cap
(C_{\vec{t}}\times C'_{\vec{t}}))\
\end{equation}
Where the union of $\vec{r},\vec{s},\vec{t}$ ranges over all
vectors such that the set of numbers appearing in their
coordinates are pairwise disjoint. This is a union of disjoint
sets. Thus:

\[\tilde{\mu}([w]_k)=\sum_{\vec{s},\vec{t},\vec{r}}(\mu_1 \times
\mu_2) ((Z\cap A_{\vec{r}} \cap B_{\vec{s}} \cap C_{\vec{t}})
\times (A'_{\vec{r}}\cap B'_{\vec{s}} \cap C'_{\vec{t}}))\]
\begin{equation} \label{mu_1_prod_mu_2}
\tilde{\mu}([w]_k)=\sum_{\vec{s},\vec{t},\vec{r}}\mu_1(Z\cap
A_{\vec{r}} \cap B_{\vec{s}} \cap
C_{\vec{t}})\mu_2(A'_{\vec{r}}\cap B'_{\vec{s}} \cap C'_{\vec{t}})
\end{equation}
 Now notice that for every $\vec{r},\vec{s},\vec{t}$ in the sum,
\[\mu_2(A'_{\vec{r}}\cap B'_{\vec{s}} \cap
C'_{\vec{t}})=m^{-(n_1+n'_2+n''_2)}=m^{-(n_1+n_2)}\] Also note
that $Z = \biguplus_{\vec{s},\vec{t},\vec{r}}(Z\cap A_{\vec{r}}
\cap B_{\vec{s}} \cap C_{\vec{t}})$, so
$\mu_1(Z)=\sum_{\vec{s},\vec{t},\vec{r}}\mu_1(Z\cap A_{\vec{r}}
\cap B_{\vec{s}} \cap C_{\vec{t}})$. Thus, equation
\ref{mu_1_prod_mu_2} can be simplified as follows:
\[\tilde{\mu}([w]_k)=\sum_{\vec{s},\vec{t},\vec{r}}\mu_1(Z\cap
A_{\vec{r}} \cap B_{\vec{s}} \cap
C_{\vec{t}})m^{-(n_1+n_2)}=\mu_1(Z)m^{-(n_1+n_2)}=(\frac{1}{2})^{|w|}m^{-(n_1+n_2)}\]
\end{proof}

\begin{theorem}\label{tilde_mu}
$\tilde{\mu}$ is a $\mathcal{T}_2$-invariant probability.
\end{theorem}
Our method of proving this is as follows:
 We define a countable set of $\mathcal{T}_2$-holonomies
 \[\mathcal{H}=\{g_{w,w',n}: \; n \in \mathbb{Z}, \; w,w' \in L(X) \; |w|=|w'|,\; w \equiv w' \mbox{(mod $M$)}, \; \}\]
 By proposition \ref{g_w_w'} bellow, we see that $\tilde{\mu}$ is
 invariant under $\mathcal{H}$. The we prove that $\mathcal{H}$
 generates $\mathcal{T}_2$, up to a $\tilde{\mu}$-null set (proposition
 \ref{H_generate_mod_mu} bellow). This will complete the proof.

\begin{lemma}\label{g_is_good}
Suppose $w,w' \in \mathcal{L}(X,n)$ with $w \equiv w'$. If $x,y
\in \Sigma^{\mathbb{Z}}$ such that $x_{[k-n,k]}=w$
$y_{[k-n,k]}=w'$ and $x_{[k-n,k]^c}=y_{[k-n,k]^c}$, then
\[x \in X \Leftrightarrow y \in X\]
\end{lemma}
\begin{proof}
Suppose $x \in X$. We have to show that for every $j
> n$ $y_{[k-j,j]} \not\equiv 0$ (mod $M$). Writing
$x_{[k-j,j]}=swt$ , we have $y_{[k-j,j]}=sw't$ and since $ w\equiv
w'$ (mod $M$), $sw't \equiv swt \not\equiv 0$ (mod $M$). This
shows $y \in X$. By replacing the roles of $y$ and $x$ we get: $y
\in X \Rightarrow x \in X$.
\end{proof}

 For $w,w' \in \mathcal{L}(X,n)$ with $w \equiv w'$ (mod $M$) and $k
\in \mathbb{Z}$. By lemma \ref{g_is_good} we can define
$g_{w,w',k}:[w]_k \mapsto [w']_k$ to be the Borel function that
changes the $n$ coordinates starting at $k$ from $w$ to $w'$.
\[g_{w,w',k}(\ldots,x_{k-1},w_0,\dots,w_{n-1},x_{k+n},\ldots)=(\ldots,x_{k-1},w'_0,\ldots,w'_{n-1},x_{k+n}\ldots)\]

\begin{proposition}\label{w_w'}
If $w \equiv w' (\mbox{mod } M)$ , $|w|=|w'|$, and $k \in
\mathbb{Z}$, then $\tilde{\mu}([w]_k)=\tilde{\mu}([w']_k)$.
\end{proposition}
\begin{proof}
By lemma \ref{tilde_mu_explicit},
$\tilde{\mu}([w]_k)=m^{-(n_1+n_2)}(\frac{1}{2})^{|w|}$. Since the
number of paired $\alpha$ in $w'$ is also $n_1$, we get that
$\tilde{\mu}([w]_k)=\tilde{\mu}([w']_k)$.
\end{proof}
\begin{proposition}\label{g_w_w'}
If $w \equiv w' (\mbox{mod } M)$, $|w|=|w'|$, and $k \in
\mathbb{Z}$, then $\tilde{\mu}$ is $g_{w,w',k}$ invariant.
\end{proposition}
\begin{proof}
First note that if $w \equiv w' (\mbox{mod } M)$ then for every
$s,t \in L(X)$ $swt \equiv sw't (\mbox{mod } M)$. This fact, along
with proposition \ref{w_w'} shows that
$\tilde{\mu}(A)=\tilde{\mu}(g_{w,w',k}(A))$ for every cylinder set
A. Since the cylinder sets generate the Borel sets, this shows
$\tilde{\mu}$ is $g_{w,w',k}$-invariant.
\end{proof}
For $w \in \mathcal{L}(X,n)$ define
\[H(w)=\sum_{i=0}^{n-1}\sum_{j=1}^{m}(\delta_{\alpha_j,w_i}-
\delta_{\beta_j,w_i})\] Also, for $x \in B_{-\infty}^{-\infty}$,
and $j >0$ define:
\[a_j(x)=\min \{ k >0 : \; H_{k+1}(x)=-j\}\]
\[b_j(x)=\max \{ k <0 : \; H_k(x)=-j\}\]
Note that for any $x \in B_{-\infty}^{-\infty}$, $(a_j(x))_{j \in
\mathbb{N}}$ is strictly increasing, and $(b_j(x))_{j \in
\mathbb{N}}$ is strictly decreasing. Also note that $x_{a_j(x)}
\in \{\beta_1,\ldots,\beta_m\}$ and $x_{b_j(x)} \in
\{\alpha_1,\ldots,\alpha_m\}$, and if $x_{a_j(x)}=\beta_i$ then
$x_{b_j(x)}=\alpha_i$.
 Let $A_c^n =\{x \in B_{-\infty}^{-\infty}:\;x_{b_j(x)}=x_{b_{j+c}(x)}\; \forall j> n\}$.

\begin{lemma}\label{mu_A_c_0}
$\tilde{\mu}(A_c^n)=0$ for all $c \in \mathbb{Z}\setminus \{0\}$,
$n \geq 0$
\end{lemma}
\begin{proof}
For $z \in S_{-\infty}^{-\infty}$ define
\[\tilde{b}_j(z)=\max \{ k < 0: \; \tilde{H}_k(z)=j\}\]
For any  $x \in B_{-\infty}^{-\infty}$,
$\tilde{b}_j(z(x))=b_j(x)$. Now, for $J \subset \mathbb{N}$ with
$|J| < \infty$:
\[ \tilde{\mu}(\{x_{b_j(x)}=x_{b_{j+c}(x)} \; \mbox{for $j\in J$
}\})=\]
\[ (\mu_2\times\mu_1)(\{(a,z) : \; a_{l_{j,1}} = a_{l_{j,2}},\;
l_{j,1}=\tilde{b}_j(z)\; l_{j,2}=\tilde{b}_{j+c}(z) \; \mbox{for
$j \in J$} \})= (\frac{1}{m})^{|J|}\]
 This follows from the definition of
$\tilde{\mu}$ as the image of a product measure, and from the fact
that $(b_j(x))_{j \in \mathbb{N}}$ is strictly monotonic, so the
$l_{j,1}$'s are all distinct, and $l_{j,1} \neq l_{j,2}$ for $j
\in J$. Thus, $\tilde{\mu}(A_c^n)=0$.
\end{proof}
\begin{proposition}\label{H_generate_mod_mu}
 There exists a double-tail invariant set $X_0 \subset X$ with $\tilde{\mu}(X_0)=1$, such the countable set
of $\mathcal{T}_2$-holonomies
\[\mathcal{H}=\{g_{w,w',n}: \; n \in \mathbb{Z}, \; w,w' \in L(X) \; |w|=|w'|,\; w \equiv w' \mbox{(mod $M$)}, \; \}\]
generates $\mathcal{T}_2(X_0)$.
\end{proposition}

\begin{proof}
Let  $X_0 = B_{-\infty}^{-\infty}\setminus \bigcup_{n,m
> 0}\bigcup_{c \neq
0}T^{-m}A_c^n$. Since $\tilde{\mu}(B_{-\infty}^{-\infty})=1$, and
$\tilde{\mu}(A_c^n)=0$ for $c \neq 0$ by the previous lemma,
$\tilde{\mu}(X_0)=1$. Also, since $B_{-\infty}^{-\infty}$ and
$\bigcup_{n,m
> 0}\bigcup_{c \neq
0}T^{-m}A_c^n$ are $\mathcal{T}_2$-invariant sets, $X_0$ is
$\mathcal{T}_2$-invariant.
We show that $\mathcal{H}$ generates $\mathcal{T}_2(X_0)$.\\
Suppose $(x,y) \in \mathcal{T}_2(X_0)$. We must show that $y=
g(x)$ for some $g \in \mathcal{H}$. $\exists n \in \mathbb{N}$ so
that $x_{[-n,n]^c}=y_{-[n,n]^c}$.
 Let $w=x_{[-n,n]}$, $w'=y_{[-n,n]}$. Let $c=H(w)-H(w')$.\\
First assume $c \neq 0$. Let $x'=T^{-n}(x)$, $y'=T^{-n}(y)$. Then
$x'_{[0,2n]^c}=y'_{[0,2n]^c}$. For all $k > 2n$,
$H_k(x')=H_k(y')+c$. Therefore, $a_j(x')=a_{j+c}(y')$ for all
$j>2n+|c|$. Also, Since $x'_{[0,2n]^c}=y'_{[0,2n]^c}$,
$H_k(x')=H_k(y')$ for all $k<0$. So $b_j(x')=b_j(y')$ for all $j
>0$.\\
 For $j>2n+|c|$, denote $x'_{a_j(x')}=\beta_i$. Then
 $x'_{b_j(x')}=\alpha_i$. Also,
 $y'_{a_{j+c}(y')}=y'_{a_j(x')}=x'_{a_j(x')}=\beta_i$, so
 $y'_{b_{j+c}(y')}=\alpha_i$. Therefore,
 $x'_{b_{j+c}(x')}=y'_{b_{j+c}(y')}=\alpha_i$. We conclude
 that $x'_{b_j(x')}=x'_{b_{j+c}(x')}$ for all $j > 2n+|c|$.
  This proves that $x \in T^{-n}A_{c}^{2n+|c|}$, but we assumed $x
\in X_0$, so this is a
contradiction, so  $c=0$.\\
Therefore, for every $k_1<-n$ and $k_2>n$, we have:
\[H_{k_1}(x)-H_{k_2}(x)= H_{k_1}(y)-H_{k_2}(y)\]
 Let $N=\min \{k \geq n : \; H_{k+1}(x)<-2n\}$, and $N'=\max \{k<-n: \;
H_k(x)=H_{N(x)+1}(x)\}$. $N$ and $N'$ are well defined for $x \in
B_{-\infty}^{-\infty}$. We have that $H_{N+1}(x)-
H_{N'}(x)=H_{N+1}(y)- H_{N'}(y)=0$, and so $x_{[N',N]} \equiv
y_{[N',N]} \equiv 0$ (mod
$M$). Thus $y=g_{x_{[N',N]},y_{[N',N]},N'}(x)$.\\
\end{proof}

\begin{proposition}\label{tilde_mu_shift}
$\tilde{\mu}$ is a shift invariant probability.
\end{proposition}
\begin{proof}
Let $[w]_k$ be a cylinder set.By lemma \ref{tilde_mu_explicit}, we
have:
\[\tilde{\mu}([w]_k)=m^{-n_1+n_2}(\frac{1}{2})^{|w|}\]
and also:
\[\tilde{\mu}(T^{-1}[w]_k)=m^{-n_1+n_2}(\frac{1}{2})^{|w|}\]
So $\tilde{\mu}(A)=\tilde{\mu}(T^{-1}[A])$ for every Borel set
$A$.
\end{proof}

\begin{proposition}\label{h_tilde_mu}
\[h_{\tilde{\mu}}(X,T)=\log(2)+\frac{1}{2}\log(m)\]
\end{proposition}
\begin{proof}
We have $h_{\tilde{\mu}}(X,T)=\lim_{n \rightarrow
\infty}h_{\tilde{\mu}}(x_0 | x_{-1},x_{-2},\ldots,x_{-n})$. Let
\[\varpi(a_1,\ldots,a_n)=\min\{H(a_1,\ldots,a_k)
: \; 0 \leq k \leq n\}\] By applying lemma
\ref{tilde_mu_explicit}, we get:
\[h_{\tilde{\mu}}(x_0 | x_{-1}=a_1,\ldots,x_{-n}=a_{n})=
\left\{ \begin{array}{ll} \log(2m) & \mbox{if
$\varpi(a_1,\ldots,a_n) \geq 0$} \\
\log(2)+\frac{1}{2}\log(m) & \mbox{if $\varpi(a_1,\ldots,a_n) <
0$}\end{array} \right.
\]
We have $h_{\tilde{\mu}}(x_0 |
x_{-1},x_{-2},\ldots,x_{-n})=\tilde{\mu}(\varpi(a_1,\ldots,a_n) <
0)(\log(2)+\frac{1}{2}\log(m))+\tilde{\mu}(\varpi(a_1,\ldots,a_n)
\geq 0)\log(2m)$. Since $\lim_{n \rightarrow
\infty}\tilde{\mu}(\varpi(a_1,\ldots,a_n) \geq 0)=0$, we have
$h_{\tilde{\mu}}(X,T)=\log(2)+\frac{1}{2}\log(m)$.
\end{proof}
By proposition \ref{h_tilde_mu},
$h_{\tilde{\mu}}(X,T)=\log(2)+\frac{1}{2}\log(m)$. Theorem
\ref{tilde_mu} together with proposition \ref{tilde_mu_shift}
provides an example of a shift invariant probability, which is
also $\mathcal{T}_2$ invariant, yet has entropy which is strictly
less than the topological entropy, for $m \geq 2$ .

\section{No other Double-Tail Invariant Probabilities}\label{no_more}
In this section we conclude that apart from the two probabilities
described in section \ref{mu_+_mu_-} and the probability defined
in section \ref{third_prob}, there are no other ergodic
double-tail
invariant probabilities for the Dyck shift.\\
By lemma \ref{mu_+_mu_-_unique} we know that there are no more
double-tail invariant probabilities  on the sets
$B^{+\infty}_{-\infty}$ and $B^{-\infty}_{+\infty}$. We also know
by the same lemma that there are no such probabilities on
$B^{\mathbb{R}}_{-\infty}$ and $B^{-\infty}_{\mathbb{R}}$.

Our next goal is to prove $\tilde{\mu}$ is unique on
$B_{-\infty}^{-\infty}$.
\begin{proposition}\label{Balance_prob_infty_infty}
Suppose $\nu$ is a $\mathcal{T}_2(B_{-\infty}^{-\infty})$
invariant probability. Then for every $w \equiv 1$ (mod $M$),
\[\nu([w]_t)= (\frac{1}{2\sqrt{m}})^{|w|}\]
\end{proposition}
\begin{proof}
Let $[w]_t$ be a balanced cylinder with $|w|=2n$. For $i<t$,
Denote:
\[M_{i,i+2N}=\{x \in X: \; x_i^{i+2N} \equiv 1 \mbox{( mod $M$)}\}\]
Since all balanced cylinders of the same length have equal $\nu$-
probability, we can calculate $\nu([w]_t \mid M_{i,i+2N})$ by
counting the number of balanced words of length $2N$, and the
number of such balanced words  with $w$ as a subword starting at
position $t-i$. The number of balanced words of length $2N$ is
\[w^m_{2N}=\frac{\left( \begin{array}{cc}2N \\ N
\end{array} \right)}{N+1}m^k\]
For a detailed calculation see  pages 69-73 of \cite{FL}, and also
lemma 4.2 in \cite{T04}. The number balanced word of length $2N$
with $w$ as a subword starting at position $t-i$ is $w^m_{2N-2n}$.
Thus,
\[\nu([w]_t \mid M_{i,i+2N})=\frac{w^m_{2N}}{w^m_{2N-2n}}\]
By an elementary calculation, we have:
\[\lim_{N \rightarrow \infty}\nu([w]_t \mid M_{i,i+2N})=
\lim_{N \rightarrow
\infty}\frac{w^m_{2N}}{w^m_{2N-2n}}=(\frac{1}{2\sqrt{m}})^n\]
Since $\nu(B_{-\infty}^{-\infty})=1$, we have
\[\nu(\bigcap_{N_0 \in \mathbb{N}}\bigcup_{i \in -\mathbb{N}}
\bigcup_{N > N_0}M_{i,i+2N})=1\] For $N_0 > n$ define a random
variable $ \chi_{N_0}(x) := \min\{ N > N_0 : \; x \in \bigcup_{i
\in -\mathbb{N}}M_{i,i+2N} \}$. We have
\[\nu([w]_t)=\sum_{N>N_0}\nu(\chi_{N_0} = N)\nu([w]_t \mid
\chi_{N_0}=N) \rightarrow(\frac{1}{2\sqrt{m}})^n\]
\end{proof}

\begin{proposition}
$\tilde{\mu}$ is the unique $\mathcal{T}_2$ invariant probability
on $B_{-\infty}^{-\infty}$.
\end{proposition}
\begin{proof}
Suppose $\nu$ is a $\mathcal{T}_2$ invariant probability on
$B_{-\infty}^{-\infty}$. By proposition
\ref{Balance_prob_infty_infty},
\begin{equation} \label{bal_cond}
 \forall w \equiv 1 \mbox{(mod $M$)} \; \nu([w])=
(\frac{1}{2\sqrt{m}})^{|w|}
\end{equation}
For $ a \in L(X)$, we say that $w \in L(X)$ is a \emph{minimal
balanced extension} of $a$, if the following conditions hold:
\begin{enumerate}
    \item There exist $l,r \in L(X)$ such that $w=lar$.
    \item $w \equiv 1$ (mod $M$)
    \item For every  $l'$ suffix of $l$ and $r'$ prefix of $r$,
    $l'ar' \equiv 1$ implies $l'ar' = w$.
\end{enumerate}
Since for every $a \in L(X)$,
 \[[a]_t =_{\nu} \biguplus\{[w]_s : \mbox{$w$ is a minimal balanced extension of $a$, with $(w_i)_{i=t-s}^{t-s+|w|}=a$}\}\]
 We have:
 \[\nu([a]_t)= \sum_{[w]_s}\nu([w]_s)=\sum_{[w]_s}\tilde{\mu}([w]_s)=\tilde{\mu}([a]_t)\]
 Where the sum ranges over minimal balanced extensions of $a$.
 This proves $\nu=\tilde{\mu}$.
 By theorem \ref{tilde_mu}, this proves $\tilde{\mu}$ is the unique double tail invariant
probability of $B_{-\infty}^{-\infty}$.
\end{proof}

Finally, we show that no invariant other double-tail invariant
probabilities exist for  Dyck.\\
Define: $\hat{p}:\Sigma^{\mathbb{Z}} \mapsto \Sigma^{\mathbb{N}}$
by $\hat{p}((x_n)_{n \in \mathbb{Z}})=(x_n)_{n \in \mathbb{N}}$.
This is a Borel mapping that maps the two-sided Dyck shift $X$
onto the
one sided Dyck shift $Y \subset \Sigma^\mathbb{N}$.\\
 Let $K_0= \{ x \in X : H_i(x) \geq 0 ,\forall i<0 \}$, and $K_i=T^{-i}(K_0))$. Notice that
$B^s_t \subset \bigcup_{i=0}^{\infty}K_i$, for $s,t \in \{\{+\infty\},\mathbb{R}\}$. \\
\begin{lemma}\label{T_2_to_T}
If $A,B \subset Y$ is are Borel sets, and $g:A \mapsto B$ is a
$\mathcal{T}(Y)$-holonomy, then there exists a
$\mathcal{T}_2(X)$-holonomy $\tilde{g}:(\hat{p}^{-1}(A)\cap K_0)
\mapsto (\hat{p}^{-1}(B)\cap K_0)$
\end{lemma}
\begin{proof}
We define $\tilde{g}:(\hat{p}^{-1}(A)\cap K_0) \mapsto
(\hat{p}^{-1}(B)\cap K_0)$ as follows:
\[\tilde{g}(x)_n = \left\{ \begin{array}{ll}
x_n & n<0 \\
g(\hat{p}(x))_n  & n\geq 0 \end{array} \right.\] We prove that
$\tilde{g}$ takes $\hat{p}^{-1}(A)\cap K_0$ into
$\hat{p}^{-1}(B)\cap K_0$. Let $x \in \hat{p}^{-1}(A)\cap K_0$.
 Since $x_n=\tilde{g}(x)_n$ for all $n <0$, we have $H_n(x)=
 H_n(\tilde{g}(x))$ for $n <0$. Because $x \in K_0$ we have
 $H_n(\tilde{g}(x)) \geq 0$ for $i < 0$. Let $y=\tilde{g}(x)$. Now we prove that  $y \in X$. Otherwise, there exist $n_1,n_2 \in \mathbb{Z}$, such that
$n_1=\min\{l<n_2:H_l(y)=H_{n_2+1}(y)\}$, and $y_{n_1}=\alpha_i$
$y_{n_2}=\beta_j$ with $i \neq j$. If $n_1,n2  <0$ then
$y_{n_1}=x_{n_1}$, $y_{n_2}=x_{n_2}$, so this contradicts the fact
that $x \in X$. If $n_1,n2  \geq 0$, then
$y_{n_1}=g(\hat{p}(x))_{n_1}$, $y_{n_2}=g(\hat{p}(x))_{n_2}$, so
this contradicts the fact that $g(\hat{p}(x)) \in Y$.\\
We remain with the case $n_1 < 0 \leq n_2$. We have $H_{n_1}(y)
\geq 0 = H_0(y)$, and $H_{n_2+1}(y)=H_{n_2}(y)-1$ (since
$y_{n_2}=\beta_j$).Also, $H_{n_2+1}(y) = H_{n_1}(y) \geq 0$. Since
$H_i(y)-H_{i+1}(y)= \pm 1$, there must be some $l >0$ such that
$H_l(y)=H_{n+1}(y)$. This contradicts the condition on $n_1,n2$.
By the definition of $\tilde{g}$,
$\hat{p}(\tilde{g}(x))=g(\hat{p}(x))$, so $\tilde{g}(x) \in
\hat{p}^{-1}(B)$. The fact that $g$ is one to one and onto
$(\hat{p}^{-1}(B)\cap K_0)$ follows from the fact that
\[\tilde{g}^{-1}(x)_n = \left\{ \begin{array}{ll}
x_n & n<0 \\
g^{-1}(\hat{p}(x))_n  & n\geq 0 \end{array} \right.\] To complete
the proof of  the lemma we must show that $(x,\tilde{g}(x)) \in
\mathcal{T}_2(X)$. Since $g$ is a $\mathcal{T}(Y)$-holonomy,
$\hat{p}(x)$ and $g(\hat{p}(x)$ only differ in a finite number of
(positive) coordinates. $x$ and $\tilde{g}(x)$ only differ in the
coordinates where $\hat{p}(x)$ and $g(\hat{p}(x))$ differ, which
is a finite set. So $(x,\tilde{g}(x)) \in \mathcal{T}_2(X)$
\end{proof}

\begin{lemma}
There are no $\mathcal{T}_2(X)$-invariant probability measures on
X supported by $B^s_t$, $s,t \in \{ \{+\infty\},\mathbb{R}\}$.
\end{lemma}
\begin{proof}
We first prove the result for $B^{\mathbb{R}}_t$,$t \in \{
\{+\infty\}, \mathbb{R}\}$. Let $K_i=T^{-i}(K_0))$. Notice that
$B^{\mathbb{R}}_t \subset
\bigcup_{i=0}^{\infty}K_i$. \\
Suppose $\mu$ is a $\mathcal{T}_2(X)$-invariant probability
supported by $B^{\mathbb{R}}_t$,where $t \in \{ \{+\infty\},
\mathbb{R}\}$, then $\mu(K_i)>0$ for some $i \ge 0$.
Without loss of generality we can assume $\mu(K_0)>0$.\\
Define a probability $\breve{\mu}$ on $Y$ by the formula:
\[\breve{\mu}(A)=\frac{\mu(\hat{p}^{-1}(A)\cap K_0)}{\mu{K_0}}\]
By lemma \ref{T_2_to_T}, $\breve{\mu}$ is a $\mathcal{T}(Y)$
invariant probability. Also, since $\mu(B^{\mathbb{R}}_t)=1$,
\[\breve{\mu}(\{y \in Y: \liminf_{n \rightarrow + \infty}H_n(y) \in  \mathbb{R}\})=1\]
 Similarly, the existence of a $\mathcal{T}_2(X)$-invariant probability supported
by $B^{+\infty}_t$,where $t \in \{ \{+\infty\} \mathbb{R}\}$ would
result in a $\mathcal{T}(Y)$-invariant probability $\breve{\mu}$
with
\[\breve{\mu}(\{y \in Y: \liminf_{n \rightarrow + \infty}H_n(y) =  +\infty \})=1\]
But in \cite{T04} it was proved that the one sided Dyck shift has
a unique $\mathcal{T}$-invariant probability, supported by
\[\{y \in Y: \liminf_{n \rightarrow + \infty}H_n(y) = -\infty\}\]

\end{proof}


\begin{thebibliography}{alpha}

\bibitem{FL}
W. Feller, An Introduction to Probabilty Theory and it's
Applications, V. 1. John Wiley \& Sons, New York, 2nd edition,
(1957).

\bibitem{PS97}
K. Petersen and K. Schmidt, Symmetric Gibbs Measures, Transactions
of the American Mathematical Society V.349 (1997) 2775-2811.

\bibitem{T04}
T. Meyerovitch, Tail Invariant Measures of the Dyck Shift,
Preprint.

\bibitem{WK74}
W. Krieger, On the Uniqueness of the Equilibruim State,
Mathematical Systems Theory 8. (1974) 97-104.



\end{thebibliography}
\end{document}